\documentstyle[12pt]{article}

\begin{document}

\begin{center}
{\Large {\bf {Jordanian twists on deformed carrier subspaces}}}
\footnote{This  work
has been partially supported by the  Russian  Foundation for Basic
Research under grants  00-01-00500 (VDL), 99-01-00101
and INTAS-99-01459 (PPK).}\\[3mm] \vskip0.5cm

{\sc Petr P. Kulish } \vskip0.25cm

{\it Steklov Mathematical Institute,}\\ {\it St.Petersburg, 191011, Russia}
\vskip0.5cm {\sc Vladimir D. Lyakhovsky} \vskip0.25cm {\it Theoretical
Department, St. Petersburg State University,}\\ {\it St. Petersburg, 198904,
Russia.} \vskip0.5cm
\end{center}

\begin{abstract}
The nontrivial subspaces with primitive coproducts are found in the deformed
universal enveloping algebras. They can form carrier spaces for
additional Jordanian twists. The latter can be used to construct sequences
of twists for algebras whose root systems contain long series of roots. The
corresponding twist for the $so(5)$ algebra is given explicitly.
\end{abstract}

\section{Introduction}

The triangular Hopf algebras ${\cal A}(m,\Delta ,S,\eta ,\epsilon ;{\cal R})$
with $R$-matrix satisfying the unitarity condition ${\cal R}_{21}{\cal R}=1$ 
, form a subclass of quasitriangular Hopf-algebras \cite{DRQG}. Many of them
can be considered as quantizations of triangular Lie bialgebras ${\bf L}$
with antisymmetric classical $r$-matrices $r=-r_{21}$ satisfying the
classical Yang-Baxter equation. The quantization is defined by a twisting
element ${\cal F}=\sum f_{\left( 1\right) }\otimes f_{\left( 2\right) }\in
{\cal A}\otimes {\cal A}$ with an expansion ${\cal F}=1+\frac 12hr+...$ \cite
{DR83}. The terms of this decomposition were defined in \cite{DR83} using
the BCH series related to the central extension of the $r$-matrix carrier
Lie algebra ${\bf L}$. In applications the knowledge of the twisting element
is highly desirable giving (twisted) ${\cal R}$-matrix ${\cal R_F}={\cal F} 
_{21}{\cal R}{\cal F}^{-1}$ and twisted coproduct $\Delta _{{\cal F}}={\cal F 
}\Delta {\cal F}^{-1}$. The explicit expressions of the twist elements $ 
{\cal F}$ were found \cite{KLM,KLO}, for the carrier algebras ${\bf L}$ with
the special properties of their triangular decompositions. The constructions
(chains of extended Jordanian twists) proposed there for the higher
dimensional carriers ${\bf L}$, were based on the effect of primitivization
of the carrier subalgebras ${\bf L}^{\prime }\subset {\bf L}$ for the
certain twists (the full canonically extended twists) performed in ${\bf L}$ 
. In this paper we present the other possibility to compose the Jordanian
and extended Jordanian twists. We show that under certain conditions when
all the coproducts in ${\bf L}$ are nontrivially twisted there exists in $U( 
{\bf L})$ the deformed primitive carrier subspace ${\bf L}_{{\cal F} 
}^{\prime }$. Thus the twist ${\cal F}$ can be composed with the next one $ 
{\cal F}^{\prime }$ defined on ${\bf L}_{{\cal F}}^{\prime }$. This effect
enlarges widely the possible applications of extended twists and chains. In
a special form the established sequences of twists can be applied to deform
the universal enveloping algebras for the classical series $B_N$ and $C_N$.
An interesting question of connections of constructed twists with recently
found quasi-Hopf twistings \cite{JIM},\cite{ARN} will be discussed in
further publications.

We hope that the constructed twisting elements for the orthogonal Lie
algebras will describe the important deformations for now actively studied
anti-de-Sitter field theories (cf e.g. \cite{FFS}) and the quantization of
coboundary bialgebra structures of conformal algebras, described by
triangular classical $r$-matrices (cf e.g. \cite{LUK}).

\section{Extended twists and the pri\-mi\-ti\-vi\-ty of deformed carrier
subspaces}

Let $g$ be a Lie algebra with the root system $\Lambda_g$, generators $ 
L_{\lambda}, \lambda \in \Lambda_g$, Cartan generators $H_{\alpha}, \alpha
\in S_{\Lambda_g}$ and the universal enveloping algebra $U(g)$. Consider the
class ${\cal B}$ of two-dimensional Borel subalgebras in $U(g)$ associated
with a fixed Cartan generator $H$ and
$$
E_{{\cal B}}=E+\xi B_iB^i
$$
where $E \equiv L_{\lambda _0}$ has the canonically normalized root $\lambda
_0$ dual to $H$, $\left[ H,E \right] =E$, and $\left\{ B_i,B^j\right\} $ is
the set of eigenvectors for $H$ with the property
\begin{equation}
\label{cond1}
\begin{array}{l}
\ \left[ H,B_i\right] =\beta _iB_i, \\
\ \left[ H,B^j\right] =\beta ^jB^j,
\end{array}
\quad \beta _i+\beta ^i=1.
\end{equation}
This guarantees that $\left[ H,E_{{\cal B}}\right] =E_{{\cal B}}$.

Let us consider the conditions sufficient for the primitivity of the
generators in ${\cal B}$. For the generator $E_{{\cal B}}$ the situation
obviously depends on the value of $\xi$:

\begin{description}
\item[I)]  $\xi =0$, in particular this is the case of the undeformed $U(g)$
where the subalgebras ${\cal B}$ are primitive;\\

\item[II)]  $\xi \neq 0$, this case can be realized only in the deformed $ 
U(g)$.
\end{description}

The primitivity of $H$ depends on the injection of the carrier subalgebras $ 
{\bf L}^{(4)}_i$ of the preceding extended Jordanian twists (see \cite{KLM})
into the carrier algebra of the final composition of twists. Consider the
Borel subalgebra $B(2) \subset {\bf L}^{(4)}$ generated by one of the
generators $L_{{\lambda _0^{\perp }}}$ with the property $\ \left[
H_{\lambda _0^{\perp}}, L_{\lambda _0^{\perp }} \right] = L_{\lambda
_0^{\perp }}$. Let $\Phi _{{\cal J}_{\perp }}$ be the Jordanian twist based
on this $B(2)$, i.e.
\begin{equation}
\label{jordfact}\Phi _{{\cal J}_{\perp }}=\exp \{H_{\lambda _0^{\perp
}}\otimes \sigma _0^{\perp }\},\quad \quad \sigma _0^{\perp }=\ln (1+L_{
\lambda _0^{\perp }}).
\end{equation}

>From now on we shall use the notation and the normalization of the
structure constants introduced in \cite{KLO}. Let $\pi _{\perp }$ be the set
of the constituent roots for ${\lambda _0^{\perp }}$:
\begin{equation}
\label{kpi}\pi _{\perp }=\left\{ \lambda ^{\prime }, \lambda ^{\prime \prime
}\,|\,\lambda ^{\prime } +\lambda ^{\prime \prime }=\lambda _0^{\perp };
\quad \lambda ^{\prime }+\lambda _0^{\perp }, \, \lambda ^{\prime \prime
}+\lambda _0^{\perp }\notin \Lambda_g\right\}
\end{equation}
For any $\lambda ^{\prime }\in \pi _{\perp }$ there must be such an element $ 
\lambda ^{\prime \prime }\in \pi _{\perp }$ that $\lambda ^{\prime }+
\lambda ^{\prime \prime }=\lambda _0^{\perp }$. So, $\pi _{\perp }$ is
naturally decomposed as
\begin{equation}
\pi _{\perp }=\pi _{\perp }^{\prime }\,\cup \,\pi _{\perp }^{\prime \prime
},\quad \quad \pi _{\perp }^{\prime }=\{\lambda ^{\prime }\},\quad \pi
_{\perp }^{\prime \prime }=\{\lambda ^{\prime \prime }\}.
\end{equation}
Consider now the sequence of extensions $\Phi _{{\cal E}_{\perp }}$ for the
Jordanian twist $\Phi _{{\cal J}_{\perp }}$,
\begin{equation}
\label{extfact}\Phi _{{\cal E}_{\perp }}=\prod_{\lambda ^{\prime }\in \pi
_{\perp }^{\prime }}\Phi _{{\cal E}_{\lambda ^{\prime }}}=\prod_{\lambda
^{\prime }\in \pi _{\perp }^{\prime }}\exp \{L_{\lambda ^{\prime }}\otimes
L_{\lambda _0^{\perp }-\lambda ^{\prime }}e^{-\frac 12\sigma _0^{\perp }}\}
\end{equation}
The complete sequence of extended twists looks like
\begin{equation}
\label{fabk}{\cal F}_{{\cal B}_{\perp }}=\Phi _{{\cal E}_{\perp }} \Phi _{ 
{\cal J}_{\perp }}
\end{equation}

In this situation the primitivity of $H$ is guaranteed by the fact that $ 
H^{*}$ is orthogonal to the only root composed by the elements of $\pi
_{\perp }$. Let us pass to the twisted coproduct for the generator $E_{{\cal  
B}}$
\begin{equation}
\label{maincop}
\begin{array}{c}
\Delta _{
{\cal F}_{{\cal B}_{\perp }}}\left( E_{{\cal B}}\right) =\Delta _{{\cal F}_{ 
{\cal B}_{\perp }}}\left( E\right) +\xi \Delta _{{\cal F}_{{\cal B}_{\perp
}}}\left( B_i\right) \Delta _{{\cal F}_{{\cal B}_{\perp }}}\left( B^i\right)
= \\ =\Delta _{\Phi _{{\cal E}_{\perp }}}\left( E\right) +\xi \Delta _{{\cal  
F}_{{\cal B}_{\perp }}}\left( B_i\right) \Delta _{{\cal F}_{{\cal B}_{\perp
}}}\left( B^i\right)
\end{array}
\end{equation}
\begin{equation}
\label{firstcop}\Delta _{\Phi _{{\cal E}_{\perp }}}\left( E\right) =\Phi _{ 
{\cal E}_{\perp }}\left( E\otimes 1+1\otimes E\right) \left( \Phi _{{\cal E}
_{\perp }}\right) ^{-1}
\end{equation}
We have the following possibilities:

\begin{description}
\item[a)]  $\lambda ^{\prime }+\lambda _0,{\lambda _0^{\perp }}-\lambda
^{\prime }+\lambda _0$ are not in $\Lambda _g$. In this case $E$ rests
primitive and this is just the case I.\\

\item[b)]  $\lambda ^{\prime }+\lambda _0,{\lambda _0^{\perp }}-\lambda
^{\prime }+\lambda _0$ are in $\Lambda _g$ but not all of them are in $\pi
_{\perp }$. In this situation the carrier subalgebra is to be enlarged. For
instance these roots may be in the other links of the chain of twists. We
shall consider this important possibility elsewhere.\\

\item[c)]  $\lambda ^{\prime }+\lambda _0,{\lambda _0^{\perp }}-\lambda
^{\prime }+\lambda _0$ are in $\pi _{\perp }$ . Let us concentrate our
attention on this case.
\end{description}

When the tensor ${\bf C}= \sum_{\lambda ^{\prime }}\left( L_{\lambda
^{\prime }}\otimes L_{{\lambda _0^{\perp }}-\lambda ^{\prime }}e^{-\frac
12\sigma _0^{\perp }}\right) $ is an invariant of the generic $B\left(
2\right) $ subalgebra (with $\left[ H,E\right] =E$ ),
$$
\ \left[ {\bf C},\Delta^{{\rm prim}}(H) \right]= \left[ {\bf C},\Delta^{{\rm  
prim}}(E) \right]= 0,
$$
we have the matreshka effect \cite{KLO} and are again in the situation I.

Let us pass to the case of a noninvariant ${\bf C}$. Moreover, we shall
suppose that the subspaces $\pi _{\perp }^{\prime }$ and $\pi _{\perp
}^{\prime \prime}$ are not conserved by the shift with the root $\lambda _0$
and impose two additional conditions: (i) $\lambda ^{\prime }+\lambda _0\in
\pi _{\perp }^{\prime \prime }$, (ii) $\lambda ^{\prime }+\lambda _0^{\perp
} $ and $\lambda ^{\prime \prime }+\lambda _0^{\perp }$ are not in $ 
\Lambda_g $ , (iii) all the $\lambda _0$-series of the roots are short. The
minimal subset involved in this action contains two roots in $\pi _{\perp
}^{\prime } $ : $\lambda ^{\prime }$ and
\begin{equation}
\label{rootcond}\widetilde{\lambda ^{\prime }}=-\lambda ^{\prime }-\lambda
_0+{\lambda _0^{\perp }.}
\end{equation}
This means that it is sufficient to regard the factors $\left\{
B_i,B^j\right\} $ as depending on the list of generators:
\begin{equation}
\label{carlist}L_{\lambda ^{\prime }},\quad L_{{\lambda _0^{\perp }}
-\lambda^{\prime }},\quad L_{\widetilde{\lambda ^{\prime }}},\quad L_{{ 
\lambda _0^{\perp }}- \widetilde{\lambda ^{\prime }}}, \quad L_{{\lambda
_0^{\perp }}}\quad {\rm and\quad } L_{{\lambda _0}}.
\end{equation}
The twisted coproducts $\Delta _{{\cal F}_{{\cal B}_{\perp }}}$ for the
first five of them are known (see \cite{KLO}). In particular, $\sigma
_0^{\perp }$ is primitive and
\begin{equation}
\label{cope1}
\begin{array}{c}
\Delta _{
{\cal F}_{{\cal B}_{\perp }}}\left( L_{\lambda ^{\prime };\widetilde{\lambda
^{\prime }}}\right) = L_{\lambda ^{\prime };\widetilde{\lambda ^{\prime }}}
\otimes e^{- \frac12 \sigma _0^{\perp }} + 1 \otimes L_{\lambda ^{\prime };
\widetilde{\lambda ^{\prime }}}, \\ \Delta _{{\cal F}_{{\cal B}_{\perp
}}}\left( L_{{\lambda _0^{\perp }} - \lambda ^{\prime };{\lambda _0^{\perp }}
- \widetilde{\lambda ^{\prime }}}\right) = L_{{\lambda _0^{\perp }}-\lambda
^{\prime };{\lambda _0^{\perp }} - \widetilde{\lambda ^{\prime }}}\otimes
e^{+ \frac12 \sigma _0^{\perp }} + e^{\sigma _0^{\perp }}\otimes L_{{\lambda
_0^{\perp }} - \lambda ^{\prime };{\lambda _0^{\perp }}-\widetilde{\lambda
^{\prime }}}.
\end{array}
\end{equation}
Returning to the expression (\ref{firstcop}) we get the last coproduct
\begin{equation}
\label{cope2}
\begin{array}{c}
\Delta _{
{\cal F}_{{\cal B}_{\perp }}}\left( E\right) =\prod_{\lambda ^{\prime },
\widetilde{\lambda ^{\prime }}}\exp \{L_{\lambda ^{\prime }}\otimes L_{{\
\lambda _0^{\perp }}-\lambda ^{\prime }}e^{-\frac 12\sigma _0^{\perp
}}\}\left( E\otimes 1+1\otimes E\right) \\ \exp \{-L_{\lambda ^{\prime
}}\otimes L_{
{\lambda _0^{\perp }}-\lambda ^{\prime }}e^{-\frac 12\sigma _0^{\perp }}\}=
\\ =E\otimes 1+1\otimes E+C_{\lambda ^{\prime }\lambda _0}^{
{\lambda _0^{\perp }}- \widetilde{\lambda ^{\prime }}}L_{{\lambda _0^{\perp } 
} -\widetilde{\lambda ^{\prime }}}\otimes L_{{\lambda _0^{\perp }}-\lambda
^{\prime }}e^{-\frac 12\sigma _0^{\perp }}+ \\ +C_{\lambda ^{\prime }\lambda
_0}^{{\lambda _0^{\perp }}- \widetilde{\lambda ^{\prime }}}L_{{\lambda
_0^{\perp }} }\otimes L_{{\lambda _0^{\perp }}- \widetilde{\lambda ^{\prime } 
}}L_{{\lambda _0^{\perp }}-\lambda ^{\prime }}e^{-\sigma _0^{\perp }} + C_{
\widetilde{ \lambda ^{\prime }} \lambda _0}^{{\lambda _0^{\perp }}-\lambda
^{\prime }}L_{{\ \lambda _0^{\perp }}-\lambda ^{\prime }}\otimes L_{{\lambda
_0^{\perp }}- \widetilde{\lambda ^{\prime }}}e^{-\frac 12\sigma _0^{\perp
}}.
\end{array}
\end{equation}
We are still free in the normalization of $E$; we fix it so that $C_{\lambda
^{\prime }\lambda _0}^{{\lambda _0^{\perp }}-\widetilde{\lambda ^{\prime }}
}=-1.$ Notice that the condition $\lambda ^{\prime }+\lambda _0={\lambda
_0^{\perp }}-\widetilde{\lambda ^{\prime }}$ means also that we get a closed
subalgebra. According to the normalization of the ${\bf L}^{(4)}$ structure
constants mentioned above this gives the following value for the last
remaining constant:
$$
\left. C_{\lambda ^{\prime }\lambda _0}^{{\lambda _0^{\perp }}- \widetilde{ 
\lambda ^{\prime }}}C_{{\lambda _0^{\perp }} -\widetilde{\lambda ^{\prime }}
\widetilde{\lambda ^{\prime }}}^{{\lambda _0^{\perp }}}+\sum_{\nu}C_{
\widetilde{\lambda ^{\prime }}\lambda ^{\prime }}^{\nu}C_{\nu\lambda _0}^{{ 
\lambda _0^{\perp }}}+C_{\lambda _0\widetilde{ \lambda ^{\prime }}}^{{ 
\lambda _0^{\perp }}-\lambda ^{\prime }}C_{{\lambda _0^{\perp }}-\lambda
^{\prime }\lambda ^{\prime }}^{{\lambda _0^{\perp }} }=0 \right\}
\Rightarrow C_{\widetilde{\lambda ^{\prime }}\lambda _0}^{{\lambda _0^{\perp
}}-\lambda ^{\prime }}=-1
$$
The final expression for the coproduct $\Delta _{{\cal F}_{{\cal B} 
_{\perp}}}\left( E\right) $ is
$$
\begin{array}{lcl}
\Delta _{{\cal F}_{{\cal B}_{\perp }}}\left( E\right) & = & E\otimes
1+1\otimes \left( E+L_{
{\lambda _0^{\perp }}-\widetilde{\lambda^{\prime }}}L_{{\lambda _0^{\perp }}
-\lambda ^{\prime }}e^{-\sigma_0^{\perp }}\right)- \\  &  & -\left( L_{
{\lambda _0^{\perp }}-\widetilde{\lambda ^{\prime }}}\otimes L_{{\lambda
_0^{\perp }} -\lambda ^{\prime }}+L_{{\lambda _0^{\perp }}- \lambda ^{\prime
}}\otimes L_{{\ \lambda _0^{\perp }}-\widetilde{\lambda ^{\prime }}}\right)
\left(1 \otimes e^{-\frac12\sigma _0^{\perp }} \right) \\  &  & -e^{\sigma
_0^{\perp }}\otimes L_{{\lambda _0^{\perp }}- \widetilde{\lambda ^{\prime }} 
}L_{{\lambda _0^{\perp }}-\lambda ^{\prime }}e^{-\sigma _0^{\perp }}.
\end{array}
$$
Thus for the factors $\left\{ B_i,B^j\right\} $ we have the following
(coalgebraic) equation:
\begin{equation}
\label{coequ}
\begin{array}{c}
1\otimes \left( L_{
{\lambda _0^{\perp }}- \widetilde{\lambda ^{\prime }}}L_{{\lambda _0^{\perp } 
} -\lambda ^{\prime }}e^{-\sigma _0^{\perp }}\right) - \\ -\left( L_{
{\lambda _0^{\perp }}-\widetilde{\lambda ^{\prime }}}\otimes L_{{\lambda
_0^{\perp }}-\lambda ^{\prime }}+L_{{\lambda _0^{\perp }}-\lambda ^{\prime
}}\otimes L_{{\lambda _0^{\perp }}-\widetilde{\lambda ^{\prime }}}\right)
\left( 1 \otimes e^{-\frac 12\sigma _0^{\perp }} \right) \\ -e^{\sigma
_0^{\perp }}\otimes L_{
{\lambda _0^{\perp }}- \widetilde{\lambda ^{\prime }}}L_{{\lambda _0^{\perp } 
} -\lambda ^{\prime }}e^{-\sigma _0^{\perp }}+\xi \Delta _{{\cal F}_{{\cal B}
_{\perp }}}\left( B_i\right) \Delta _{{\cal F}_{{\cal B}_{\perp }}}\left(
B^i\right) = \\ \left( \xi B_iB^i\right) \otimes 1+1\otimes \left( \xi
B_iB^i\right) .
\end{array}
\end{equation}
It demonstrates that the list of realizations for $B$'s in (\ref{carlist})
might be reduced to the set
$$
L_{{\lambda _0^{\perp }}-\widetilde{\lambda ^{\prime }}}, \quad L_{{ 
\lambda_0^{\perp }}-\lambda ^{\prime }}\quad {\rm and\quad } L_{{\lambda
_0^{\perp }}}.
$$
This immediately leads to the solution of the equations (\ref{cope1}) and ( 
\ref{coequ}) that gives the following answer:
\begin{equation}
\label{newbor}E_{{\cal B}}=E+L_{{\lambda _0^{\perp }}-\widetilde{\lambda
^{\prime }}}L_{{\lambda _0^{\perp }}- \lambda ^{\prime }}e^{-\sigma
_0^{\perp}}; \qquad \xi=1.
\end{equation}

We have proved that in the twisted universal enveloping algebra $U_{{\cal F} 
_{{\cal B}_{\perp }}}$ one can find the deformed primitive Borel subspace.
In the case I this was one of the basic effects in constructing chains of
extended twists \cite{KLO}. In the case II this also provides the
possibility to compose new twists such as
\begin{equation}
\label{newtwist}{\cal F}_{{\cal BJE}}=\Phi _{{\cal BJ}}\Phi _{{\cal E}
_{\perp }}\Phi _{{\cal J}_{\perp }}
\end{equation}
where the second Jordanian factor is defined on the deformed carrier
subspace generated by $\left\{ H, E_{{\cal B}} \right\}$,
\begin{equation}
\label{newjord}
\begin{array}{c}
\Phi _{
{\cal BJ}}=\exp \left( H\otimes \sigma \left( E_{{\cal B}}\right) \right) ,
\\ \sigma \left( E_{{\cal B}}\right) =\ln \left( 1+E_{{\cal B}}\right)
\end{array}
\end{equation}
The carrier algebra for the twists ${\cal F}_{{\cal BJE}}$ is 8-dimensional
with the generators
\begin{equation}
\label{8-list} \left\{ H, \quad H_{\lambda _0^{\perp }}, \quad L_{{\lambda _0 
}}, \quad L_{{\lambda_0^{\perp }}}, \quad L_{\lambda ^{\prime }}, \quad L_{{ 
\lambda _0^{\perp }}-\lambda^{\prime }}, \quad L_{\widetilde{\lambda
^{\prime }}}, \quad L_{{\lambda _0^{\perp }}- \widetilde{\lambda ^{\prime }} 
}\right\}
\end{equation}
and the set of roots
\begin{equation}
\label{6-root} \lambda _0, \quad \lambda _0^{\perp }, \quad \lambda ^{\prime
}, \quad {\lambda _0^{\perp }}-\lambda ^{\prime }, \quad \widetilde{\lambda
^{\prime }}, \quad {\lambda_0^{\perp }}-\widetilde{\lambda ^{\prime }}.
\end{equation}

It is necessary to distinguish the following special case of the general
structure considered above. This is the case when in the root condition (\ref
{rootcond}) the roots $\lambda ^{\prime }$ and $\widetilde{\lambda ^{\prime } 
}$ coincide:
\begin{equation}
\label{coinc} \lambda ^{\prime }+\lambda _0={\lambda _0^{\perp }-}\lambda
^{\prime }
\end{equation}
The root subset now contains 4 roots:
\begin{equation}
\label{4-root} \lambda_0^{\perp },\quad \lambda_0,\quad \lambda^{\prime },
\quad \lambda_0^{\perp }-\lambda^{\prime }
\end{equation}
and we have a long $\lambda ^{\prime }$-series of the root ${\lambda
_0^{\perp }:} $
\begin{equation}
\lambda_0^{\perp }, \quad \lambda_0^{\perp }- \lambda^{\prime }, \quad
\lambda_0^{\perp }-2\lambda^{\prime } \quad \in \Lambda_g
\end{equation}
This is obviously the property characteristic for the series $B_N$ and $C_N$
of simple Lie algebras and the exceptional algebra $F_4$. In such
deformations the twisting elements are still of the form (\ref{newtwist})
and (\ref{newjord}) while the expression (\ref{newbor}) for $E_{{\cal B}}$
is to be substituted by
$$
E_{{\cal BO}}=E+ \frac12 \left( L_{{\lambda _0^{\perp }} -\lambda ^{\prime
}}\right)^2e^{-\sigma _0^{\perp }}.
$$

\section{Examples}

Identifying the roots (\ref{6-root}) or (\ref{4-root}) of the carrier
subalgebras with the root subsets in simple Lie algebras we can perform new
twisting deformations for them specific for the properties of their root
systems.

In the case of 8-dimansional algebra ${\bf L}$ the minimal simple algebra
where the effect described above can be illustrated is the Lie algebra $ 
U(sl(4))$. The elements (\ref{8-list}) can be identified with the following
generators of $sl(4)$:
$$
\begin{array}{c}
H=H_{23},\quad H_{\lambda _0^{\perp }}=H_{14},\quad L_{
{\lambda _0}}=E_{23},\quad L_{{\lambda _0^{\perp }}}=E_{14}, \\ L_{\lambda
^{\prime }}=E_{12},\quad L_{{\lambda _0^{\perp }}-\lambda ^{\prime
}}=E_{24},\quad L_{\widetilde{\lambda ^{\prime }}}=-E_{34},\quad L_{{\lambda
_0^{\perp }}-\widetilde{\lambda ^{\prime }}}=E_{13}.
\end{array}
$$

After the first (minimal) extended twist \cite{KLM}
$$
\Phi _{{\cal E}}\Phi _{{\cal J}_{\perp }}=e^{\left( E_{12}\otimes
E_{24}e^{-\frac 12\sigma _{14}}\right) }e^{\left( H_{14}\otimes \sigma
_{14}\right) },\qquad \sigma _{14}=\ln \left( 1+E_{14}\right)
$$
the second extension factor can have the form
$$
\Phi _{{\cal E^{\prime }}}=\exp \left( -E_{34}\otimes E_{13}e^{-\frac
12\sigma _{14}}\right)
$$
According to the formula (\ref{newbor}) after the sequence of extended
twists $\Phi _{{\cal E^{\prime }}}\Phi _{{\cal E}}\Phi _{{\cal J}_{\perp }}$
the deformed carrier subspace of primitive elements forms the Borel
subalgebra ${\cal B}=\left\{ H_{23},E_{{\cal B}}\right\} \ E_{{\cal B} 
}\equiv -E_{23}+E_{13}E_{24}e^{-\sigma _{14}}$. So it is possible to use $ 
\sigma _{{\cal B}}=\ln \left( 1+E_{{\cal B}}\right) $ and apply additionally
the twist ( see (\ref{newjord})) $\Phi _{{\cal BJ}}=\exp \left(
H_{23}\otimes \sigma _{{\cal B}}\right) $ to the deformed algebra $U_{{\cal  
E^{\prime }EJ}}(sl(4))$. In this case the final twisting element looks like
\begin{equation}
\label{sl4jeej}{\cal F}_{{\cal BJE^{\prime }EJ}}=\Phi _{{\cal BJ}}\Phi _{ 
{\cal E^{\prime }}}\Phi _{{\cal E}}\Phi _{{\cal J}_{\perp }}.
\end{equation}
We shall consider such cases in full details in the forthcoming publications.

The minimal simple Lie algebra that have the root subset (\ref{4-root}) is $ 
so(5)$. When the root system of $so(2N+1)$ is fixed in the standard $e$ 
-basis as
$$
\Lambda _{so(2N+1)}=\left\{ \pm e_i, \; \pm e_i\pm e_j \,|i,j=1,\ldots
,N\right\}
$$
then in accordance with the property (\ref{coinc}) the set (\ref{4-root})
can be injected into $\Lambda _{so(5)}$ as follows
$$
\lambda _0^{\perp }=e_1+e_2,\quad \lambda _0=e_1-e_2, \quad \lambda ^{\prime
}=e_2,\quad {\lambda _0^{\perp }}-\lambda ^{\prime }=e_1.
$$
Thus we get the 6-dimensional subalgebra ${\bf L}^{\left( 6\right) }\subset
so(5)$ generated by the set
$$
\left\{ H_{1+2},\,E_{1+2},\,H_{1-2},\, E_{1-2},\,E_1,\,E_2\right\} .
$$
In terms of the ordinary antisymmetric Okubo matrices $M_{ij}$ the following
list of generators in the defining representation $d\left( {\bf L}^{\left( 6
\right) }\right) $,
$$
\begin{array}{c}
d\left( H_{\lambda _0^{\perp }}\right) =d\left( H_{1+2}\right) =-\frac
i2\left( M_{12}+M_{34}\right) , \\
d\left( H\right) =d\left( H_{1-2}\right) =-\frac i2\left(
M_{12}-M_{34}\right) , \\
d\left( L_{
{\lambda _0^{\perp }}}\right) =d\left( E_{1+2}\right) =\frac 12\left(
-M_{24}+iM_{23}+iM_{14}+M_{13}\right) , \\ d\left( L_{
{\lambda _0}}\right) =d\left( E_{1-2}\right) =\frac 12\left(
-M_{24}-iM_{23}+iM_{14}-M_{13}\right) , \\ d\left( L_{
{\lambda _0^{\perp }}-\lambda ^{\prime }}\right) =d\left( E_1\right) =\frac
1{\sqrt{2}}\left( M_{25}-iM_{15}\right) , \\ d\left( L_{\lambda ^{\prime
}}\right) =d\left( -E_2\right) =\frac 1{\sqrt{2}}\left(
-M_{45}+iM_{35}\right) ,
\end{array}
$$
fits the normalization conditions for ${\bf L}^{\left( 6\right) }$.

The canonical extended twist ${\cal F_{EJ}}$ based on ${\bf L}^{\left(
4\right) }$ with the generators $\left\{ H_{1+2},E_{1+2},E_1,-E_2\right\} $,
$$
{\cal F_{EJ}}=\exp \left( -E_2\otimes E_1e^{-\frac 12\sigma _{1+2}}\right)
\exp \left( H_{1+2}\otimes \sigma _{1+2}\right) ,\quad \sigma _{1+2}=\ln
\left( 1+E_{1+2}\right) ,
$$
leads to the deformed algebra $U_{{\cal EJ}}\left( {\bf L}^{\left( 6\right)
}\right) $ with the coproducts:
$$
\begin{array}{c}
\Delta _{
{\cal EJ}}\left( H_{1+2}\right) =H_{1+2}\otimes e^{-\sigma _{1+2}}+1\otimes
H_{1+2}+E_2\otimes E_1e^{-\frac 32\sigma _{1+2}}, \\ \Delta _{
{\cal EJ}}\left( E_{1+2}\right) =E_{1+2}\otimes e^{\sigma _{1+2}}+1\otimes
E_{1+2}, \\ \Delta _{
{\cal EJ}}\left( E_2\right) =E_2\otimes e^{-\frac 12\sigma _{1+2}}+1\otimes
E_2, \\ \Delta _{
{\cal EJ}}\left( E_1\right) =E_1\otimes e^{\frac 12\sigma _{1+2}}+e^{\sigma
_{1+2}}\otimes E_1, \\ \Delta _{
{\cal EJ}}\left( H_{1-2}\right) =H_{1-2}\otimes 1+1\otimes H_{1-2}, \\
\Delta _{
{\cal EJ}}\left( E_{1-2}\right) =E_{1-2}\otimes 1+1\otimes
E_{1-2}-E_1\otimes E_1e^{-\frac 12\sigma _{1+2}}- \\ -\frac 12E_{1+2}\otimes
E_1^2e^{-\sigma _{1+2}}.
\end{array}
$$
According to the arguments presented in Sect.2 we have in $U_{{\cal EJ} 
}\left( {\bf L}^{\left( 6\right) }\right) $ the primitive subalgebra ${\cal B 
}=\left\{ H_{1-2},E_{{\cal BO}}\right\} $ on the deformed subspace with
$$
E_{{\cal BO}}=E_{1-2}+\frac 12E_1^2e^{-\sigma _{1+2}}.
$$
In this case the ''shifted'' Jordanian factor (see (\ref{newjord}))
$$
\begin{array}{c}
\Phi _{
{\cal BJ}}=\exp \left( H\otimes \sigma _{{\cal BO}}\right) =\exp \left(
H_{1-2}\otimes \sigma _{{\cal BO}}\right) , \\ \sigma _{{\cal BO}}=\ln
\left( 1+E_{{\cal BO}}\right) ,
\end{array}
$$
can be applied to $U_{{\cal EJ}}\left( {\bf L}^{\left( 6\right) }\right) $
and/or to $U_{{\cal EJ}}\left( so(5)\right) $. The result will be the
twisted $U_{{\cal BJEJ}}\left( so(5)\right) \supset U_{{\cal BJEJ}}\left(
{\bf L}^{\left( 6\right) }\right) $ with the costructure defined by the
relations:
$$
\begin{array}{lcl}
\Delta _{{\cal BJEJ}}\left( H_{1+2}\right)  & = &
\begin{array}{l}
H_{1+2}\otimes e^{-\sigma _{1+2}}+1\otimes H_{1+2}+ \\
+E_2\otimes E_1e^{-\frac 32\sigma _{1+2}-\frac 12\sigma _{{\cal BO}}}-\frac
12H_{1-2}\otimes E_1^2e^{-2\sigma _{1+2}-\sigma _{{\cal BO}}},
\end{array}
\\
\Delta _{{\cal BJEJ}}\left( H_{1-2}\right)  & = & H_{1-2}\otimes e^{-\sigma
_{
{\cal BO}}}+1\otimes H_{1-2}, \\ \Delta _{{\cal BJEJ}}\left( E_{1+2}\right)
& = & E_{1+2}\otimes e^{\sigma _{1+2}}+1\otimes E_{1+2}, \\
\Delta _{{\cal BJEJ}}\left( E_2\right)  & = & E_2\otimes e^{-\frac 12\sigma
_{1+2}-\frac 12\sigma _{
{\cal BO}}}+1\otimes E_2-H_{1-2}\otimes E_1e^{-\sigma _{1+2}-\sigma _{{\cal  
BO}}}, \\ \Delta _{{\cal BJEJ}}\left( E_1\right)  & = & E_1\otimes e^{\frac
12\sigma _{1+2}+\frac 12\sigma _{
{\cal BO}}}+e^{\sigma _{1+2}}\otimes E_1, \\ \Delta _{{\cal BJEJ}}\left(
E_{1-2}\right)  & = &
\begin{array}{l}
E_{1-2}\otimes e^{\sigma _{
{\cal BO}}}+1\otimes E_{1-2}- \\ -E_1\otimes E_1e^{-\frac 12\sigma
_{1+2}+\frac 12\sigma _{{\cal BO}}}-\frac 12E_{1+2}\otimes E_1^2e^{-\sigma
_{1+2}},
\end{array}
\\
\Delta _{{\cal BJEJ}}\left( E_{{\cal BO}}\right)  & = & E_{{\cal BO}}\otimes
e^{\sigma _{{\cal BO}}}+1\otimes E_{{\cal BO}},
\end{array}
$$
$$
\begin{array}{l}
\Delta _{
{\cal BJEJ}}\left( E_{-1}\right) =E_{-1}\otimes e^{-\frac 12\left( \sigma _{ 
{\cal BO}}+\sigma _{1+2}\right) }+1\otimes E_{-1} \\ \qquad +H_{1-2}\otimes
\left( +E_{-2}+2H_1E_1e^{-\sigma _{1+2}}+\frac 12E_1^2E_2e^{-2\sigma
_{1+2}}\right) e^{-\sigma _{
{\cal BO}}} \\ \qquad +\frac 12H_{1-2}\otimes E_1\left( \left(
E_1^2e^{-2\sigma _{1+2}}-2\right) e^{-\sigma _{
{\cal BO}}}+1\right) e^{-\left( \sigma _{{\cal BO}}+\sigma _{1+2}\right) }+
\\ \qquad -\frac 12H_{1-2}^2
{\bf \otimes }E_1\left( \left( E_1^2e^{-2\sigma _{1+2}}-2\right) e^{-\sigma
_{{\cal BO}}}+2\right) e^{-\left( \sigma _{{\cal BO}}+\sigma _{1+2}\right) }
\\ \qquad +H_{1-2}E_2\otimes \left( \left( \frac 32E_1^2e^{-2\sigma
_{1+2}}-1\right) e^{-\sigma _{
{\cal BO}}}+1\right) e^{-\frac 12\left( \sigma _{{\cal BO}}+\sigma
_{1+2}\right) } \\ \qquad -H_{1+2}\otimes E_2e^{-\sigma
_{1+2}}+H_{1+2}H_{1-2}\otimes E_1e^{-\left( \sigma _{
{\cal BO}}+2\sigma _{1+2}\right) } \\ \qquad +E_{2-1}\otimes E_1e^{-\left(
\sigma _{
{\cal BO}}+\sigma _{1+2}\right) } \\ \qquad +H_{1+2}E_2\otimes \left(
1-e^{-\sigma _{1+2}}\right) e^{-\frac 12\left( \sigma _{
{\cal BO}}+\sigma _{1+2}\right) } \\ \qquad -E_2\otimes \left(
2H_1+E_1E_2e^{-\sigma _{1+2}}\right) e^{-\frac 12\left( \sigma _{
{\cal BO}}+\sigma _{1+2}\right) } \\ \qquad +E_2^2\otimes E_1\left( \frac
12-e^{-\sigma _{1+2}}\right) e^{-\left( \sigma _{{\cal BO}}+\sigma
_{1+2}\right) },
\end{array}
$$
$$
\begin{array}{l}
\Delta _{
{\cal BJEJ}}\left( E_{-2}\right) =E_{-2}\otimes e^{+\frac 12\left( \sigma _{ 
{\cal BO}}-\sigma _{1+2}\right) }+1\otimes E_{-2} \\ \qquad +E_2\otimes
\left( 1-e^{-\sigma _{{\cal BO}}}\right) e^{+\frac 12\left( \sigma _{{\cal BO 
}}-\sigma _{1+2}\right) }+H_{1-2}\otimes E_1e^{-\left( \sigma _{{\cal BO} 
}+\sigma _{1+2}\right) },
\end{array}
$$
$$
\begin{array}{l}
\Delta _{
{\cal BJEJ}}\left( E_{2-1}\right) =E_{2-1}\otimes e^{-\sigma _{{\cal BO} 
}}+1\otimes E_{2-1} \\ \qquad +H_{1-2}\otimes \left\{
\begin{array}{c}
\frac 12\left( 1+e^{-\sigma _{1+2}}\right) -e^{-\sigma _{
{\cal BO}}} +2H_{1-2}+
\\ +E_1e^{-\sigma _{1+2}}\left( \frac 12E_1e^{-\left(
\sigma _{{\cal BO}}+\sigma _{1+2}\right) }+E_2\right)
\end{array}
\right\} e^{-\sigma _{
{\cal BO}}} \\ \qquad +H_{1-2}E_2\otimes E_1e^{-\frac 32\left( \sigma _{
{\cal BO}}+\sigma _{1+2}\right) } \\ \qquad +H_{1-2}^2\otimes \left( \left(
1-\frac 12E_1^2e^{-2\sigma _{1+2}}\right) e^{-\sigma _{
{\cal BO}}}-1\right) e^{-\sigma _{{\cal BO}}} \\ \qquad -E_2\otimes
E_2e^{-\frac 12\left( \sigma _{{\cal BO}}+\sigma _{1+2}\right) }+\frac
12E_2^2\otimes \left( 1-e^{-\sigma _{1+2}}\right) e^{-\sigma _{{\cal BO}}},
\end{array}
$$

$$
\begin{array}{l}
\Delta _{
{\cal EBJJ}}\left( E_{-1-2}\right) =E_{-1-2}\otimes e^{-\sigma
_{1+2}}+1\otimes E_{-1-2}+ \\ \qquad +H_{1-2}\otimes \left( \frac 12\left(
e^{-\sigma _{
{\cal BO}}}-1\right) -E_1E_{-2}e^{-\sigma _{{\cal BO}}}\right) e^{-\sigma
_{1+2}} \\ \qquad +H_{1-2}\otimes \left\{
\begin{array}{c}
e^{-\sigma _{1+2}}+\frac 12e^{-\sigma _{
{\cal BO}}}-H_{1+2}- \\ -\frac 14\left( E_1^2e^{-\left( \sigma _{{\cal BO} 
}+2\sigma _{1+2}\right) }+1\right)
\end{array}
\right\} E_1^2e^{-\left( \sigma _{
{\cal BO}}+2\sigma _{1+2}\right) } \\ \qquad +\frac 12H_{1-2}^2\otimes
\left( 1-e^{-\sigma _{
{\cal BO}}}+\frac 12E_1^2e^{-\left( \sigma _{{\cal BO}}+2\sigma
_{1+2}\right) }\right) E_1^2e^{-\left( \sigma _{{\cal BO}}+2\sigma
_{1+2}\right) } \\ \qquad +\left( H_{1+2}^2-H_{1+2}\right) \otimes \left(
e^{-\sigma _{1+2}}-1\right) e^{-\sigma _{1+2}}+2H_{1+2}\otimes
H_{1+2}e^{-\sigma _{1+2}} \\
\qquad -H_{1+2}H_{1-2}\otimes E_1^2e^{-\left( \sigma _{
{\cal BO}}+3\sigma _{1+2}\right) } \\ \qquad -E_{-1}\otimes E_1e^{-\frac
12\left( \sigma _{
{\cal BO}}+3\sigma _{1+2}\right) }-\frac 12E_{2-1}\otimes E_1^2e^{-\left(
\sigma _{{\cal BO}}+2\sigma _{1+2}\right) } \\ \qquad +E_2\otimes \left(
E_{-2}+\left( 2H_{1+2}-2e^{-\sigma _{1+2}}+1\right) E_1e^{-\sigma
_{1+2}}\right) e^{-\frac 12\left( \sigma _{
{\cal BO}}+\sigma _{1+2}\right) } \\ \qquad +E_2^2\otimes \left( \frac
12e^{-\sigma _{1+2}}\left( 1-e^{-\sigma _{
{\cal BO}}}\right) +E_1^2e^{-\left( \sigma _{{\cal BO}}+2\sigma
_{1+2}\right) }\left( e^{-\sigma _{1+2}}-\frac 14\right) \right)  \\ \qquad
+H_{1-2}E_2\otimes \left( e^{-\sigma _{
{\cal BO}}}-1-E_1^2e^{-\left( \sigma _{{\cal BO}}+2\sigma _{1+2}\right)
}\right) E_1e^{-\frac 12\left( \sigma _{{\cal BO}}+3\sigma _{1+2}\right) }
\\ \qquad +H_{1+2}E_2\otimes \left( 2e^{-\sigma _{1+2}}-1\right)
E_1e^{-\frac 12\left( \sigma _{{\cal BO}}+3\sigma _{1+2}\right) }.
\end{array}
$$
Using a particular set of generators nonlinearly related with the undeformed
ones, it can be demonstrated that the bialgebra structure of $U_{{\cal BJEJ} 
}\left( so(5)\right) $ coincides with the one presented in \cite{HER} where
it was obtained as a direct solution of the conditions of coassociativity
while the twisted character of the deformation was not studied.

The Jordanian twists on the deformed carrier spaces found above, will give
rise to different new constructions enlarging the list of explicit solutions
of the Yang-Baxter equation, deformed Yangians and integrable models. 
In particular, due to the embedding of the simple 
Lie algebras $g$ into the corresponding Yangians
(as Hopf subalgebras) $U(g)\subset {\cal Y}(g)$ \cite{DRQG} the Yangian $R$
-matrix can be twisted by the same ${\cal F}$ defined for $g$. As a result 
for the case of orthogonal algebra $g=so(M)$ the $R$-matrix of  ${\cal Y}(g)$
(in the defining representation 
$\rho  \subset {\rm Mat}(M,{\bf C})\otimes {\rm Mat}(M,{\bf C})$ ) 
can be changed:
$$
\begin{array}{c}
u\rho \left( 1\otimes 1\right) +
{\cal P}-\frac u{u-1+M/2}{\cal K\quad \longrightarrow } \\ 
\quad u\rho \left( {\cal F}_{21}{\cal F}^{-1}\right) +{\cal P}-\frac
u{u-1+M/2}\rho \left( {\cal F}_{21}\right) {\cal K}\rho \left( {\cal F} 
^{-1}\right)
\end{array}. 
$$
(Here $u$ is a spectral parameter and the operator ${\cal K}$ that is
obtained from ${\cal P}$ by transposing its first tensor factor.) Henceforth
the density of the integrable spin chain hamiltonian is changing as well 
${\cal P}\rho
\left( {\cal F}_{21}{\cal F}^{-1}\right) + \frac
1{1-M/2}\rho \left( {\cal F}\right) {\cal K}\rho 
\left( {\cal F}^{-1}\right)$ 
(cf. the $sl(2)$-case \cite{KS}).

\end{document}